\documentclass[letterpaper, 10 pt, conference]{ieeeconf}  

\IEEEoverridecommandlockouts                             

\usepackage{flushend}
\usepackage{tikz}

\usepackage{amsthm}
\usepackage{amsmath}
\usepackage{amssymb}

\usepackage[inkscapearea=page]{svg}
\usepackage{adjustbox}

\usepackage{graphicx}    

\newcommand\copyrighttext{%
  \footnotesize \textcopyright 2023 IEEE. Personal use of this material is permitted.
  Permission from IEEE must be obtained for all other uses, in any current or future
  media, including reprinting/republishing this material for advertising or promotional
  purposes, creating new collective works, for resale or redistribution to servers or
  lists, or reuse of any copyrighted component of this work in other works.}
\newcommand\copyrightnotice{%
\begin{tikzpicture}[remember picture,overlay]
\node[anchor=south,yshift=10pt] at (current page.south) {\fbox{\parbox{\dimexpr\textwidth-\fboxsep-\fboxrule\relax}{\copyrighttext}}};
\end{tikzpicture}%
}

\DeclareMathOperator{\Ima}{Im}

\newtheorem{theorem}{Theorem}
\newtheorem{lemma}[theorem]{Lemma}
\newtheorem{remark}{Remark}
\newtheorem{proposition}[theorem]{Proposition}

\newtheorem{assumption}{Assumption}

\title{\LARGE \bf
Semi-Global Practical Extremum Seeking with Practical Safety*
}

\author{Alan Williams$^{1,2}$, Miroslav Krstic$^{1}$ and Alexander Scheinker$^{2}$
\thanks{*This work is supported by Los Alamos National Lab LDRD DR
Project 20220074DR}
\thanks{$^{1}$Alan Williams and Miroslav Krsti\'{c} are with the Department of Mechanical and Aerospace Engineering, University of California, San Diego, CA 92093-0411, USA, {\tt\small \{awilliam,krstic\}@ucsd.edu}.}
\thanks{$^{2}$Alan Williams and Alexander Scheinker are with Los Alamos National Lab, Los Alamos, NM 87545, USA {\tt\small ascheink@lanl.gov}}
}

\begin{document}

\maketitle
\copyrightnotice
\thispagestyle{empty}
\pagestyle{empty}

\begin{abstract}
We introduce a type of safe extremum seeking (ES) controller, which minimizes an unknown objective function while also maintaining practical positivity of an unknown barrier function. We show semi-global practical asymptotic stability of our algorithm and present an analogous notion of practical safety. The dynamics of the controller are inspired by the quadratic program (QP) based safety filter designs which, in the literature, are more commonly used in cases where the barrier function is known. Conditions on the barrier and objective function are explored showing that non convex problems can be solved. A Lyapunov argument is proposed to achieve the main results of the paper. Finally, an example is given of the algorithm which solves the constrained optimization problem.
\end{abstract}

\section{Introduction}
This paper presents an ES algorithm which can be used to minimize an unknown objective function $J(\theta)$ over parameters $\theta$ while also keeping the system safe. Safety is considered to be a measured unknown function $h(\theta)$ and is maintained by keeping $h$ positive. The analysis is based on the framework in \cite{nesic2010unifying}, preceded by notable papers \cite{tan2005non}, \cite{tan2006non} proving semi-global practical asymptotic (SPA) stability properties of extremum seeking. We use the basic ideas of \cite{nesic2010unifying} but instead consider the constrained optimization problem, and present a Lyapunov function for showing SPA stability of the reduced constrained dynamics. Additionally we present a notion analogous to that of practical stability, called `practical safety'. The dynamics we present are based on the QP safety filter \cite{ames2016control} and therefore it gives the designer a choice to weight the importance of the objective versus the safety, through a parameter $c$, while the trajectory of $\theta$ tracks toward the constrained optimum. Our algorithm approximately solves the following problem:
\begin{equation}
\min_{\theta(t)}  J(\theta(t)) \text{ s.t. } h(\theta(t)) \geq 0 \text{ for all } t \in [0, \infty) . \end{equation}

Solving constrained optimization problems using ES has been explored in other contexts in the literature, both using the framework described in \cite{nesic2010unifying} as well as other methods such as \cite{guay2015constrained}. The results in \cite{poveda2015shahshahani} prove SPA stability of an ES based controller which not only converges to set, but also an optimum constrained to a set. Other dynamics have also been considered in similar settings such as \cite{labar2019constrained} and \cite{liao2019constrained}. Switching ES algorithms have also been used to solve constrained optimization problems \cite{chen2023continuoustime}. 

Safe optimization of systems is relevant in many areas. Particle accelerators are complex, time varying systems which must be constantly tuned for optimal operation. Because of its model-independence, robustness to noise, and ability to handle large numbers of parameters simultaneously, ES is especially well suited for particle accelerator applications. For example, in \cite{scheinker2020online}, the authors used a bounded form of ES \cite{scheinker2014extremum} for real-time multi-objective optimization of a particle accelerator beam in which analytic bounds are guaranteed on each parameter update despite acting on an analytically unknown and noisy measurement function. In \cite{scheinker2020online} the algorithm was pushed towards achieving safety by weighing a measure of how far the beam was off from a prescribed optimal trajectory. Such a cost-weighing ES approach requires hand tuning of weights associated with safety and offers no guarantees that safety will be maintained.

Authors in \cite{kirschner2022tuning} tune the several beamlines using Bayesian optimization and do so without violating safety constraints like beam loss. Here ``safety'' means that the particle beam does not damaging key accelerator components during tuning based on hand-picked safety margins which result in a trade-off between safety and performance. Authors in \cite{sui2015safe} apply constrained optimization schemes using Gaussian processes to recommender systems and therapeutic spinal cord stimulation. In these examples, it is desired that recommendations to users (for movies) must not be heavily disliked, and stimulation patterns for patients must not exceed a certain pain threshold.

This work is heavily inspired by recent work \cite{williams2022practically}. Here, the QP based modification of standard ES was introduced, which approximately solves the constrained optimization problem. It was shown that from the QP safety formulation, an additive safety term can be introduced in the parameter's dynamics. The algorithm in \cite{williams2022practically} was shown to be a practically safe scheme, but only locally and under restrictive assumptions on $J$ and $h$. In this work, we explore more general assumptions on $J$ and $h$, in n-dimensions, and present a semi-global result for both safety and stability. The dynamics presented in this work are slightly different from that of \cite{williams2022practically}, due to the time scaling by $k \omega_f$, and the parameter $M^+$. The parameter $M^+$ is chosen large to bound the dynamics if the estimate of $\nabla h$ goes to zero. This is because the safety term in its exact form contains $\nabla h^T \nabla h$ in the denominator. The gain $k \omega_f$ scales the dynamics of the optimization parameter $\hat \theta$. This provides a tuning knob which has can adjust the timescale of $\hat \theta$ proportionally, a requirement for SPA stability, which allows the estimator states to converge arbitrarily quickly. Additionally, our Lyapunov analysis provides a stability result which includes non convex $h$ and $J$, although $J$ having a unique minimum on the safe set is required for SPA stability results. We provide a 2D example with a non convex $h$ to conclude. 

\textbf{Notation:} 
For a differentiable function $Q: \mathbb{R}^n \to \mathbb{R}$ we denote the gradient $\nabla Q: \mathbb{R}^n \to \mathbb{R}^n$ as the vector $\nabla Q(x) = [\partial Q(x) / \partial x_1, \partial Q(x) / \partial x_2, ..., \partial Q(x) / \partial x_n]^T$ and where the ith component is $\nabla_i Q(x) = \partial Q(x) / \partial x_i$. For $v \in \mathbb{R}^n$, the notation $|| v ||$ denotes the Euclidean norm. 
The continuous function $\beta: \mathbb{R}_{\geq 0} \to \mathbb{R}_{\geq 0}$ is of class $\mathcal{K}$ if $\beta(0)=0$ and it is strictly increasing. The continuous function $\beta: \mathbb{R}_{\geq 0} \times \mathbb{R}_{\geq 0} \to \mathbb{R}_{\geq 0}$ is of class $\mathcal{KL}$ if it is strictly
increasing in its first argument and strictly decreasing to zero in its second argument. The image of a function $h$ is denoted by $\Ima(h)$. The compact ball around a point $p$ is $B_r(p) = \{ \theta \in \mathbb{R}^n : \| \theta - p \| \leq r \} $. 
We use the term ``SPA stability'' to refer to the notion of semi-global practical asymptotic stability \cite{tan2005non}. A function $f(x,\epsilon)$ is $O(\epsilon)$ if for any compact set $\Omega$ there exists a positive pair $(\epsilon^*, k)$ such that $||f(x,\epsilon)||\leq k \epsilon$ for all $\epsilon \in (0,\epsilon^*]$ for all $x \in \Omega$.

\section{Algorithm Design}
For a better understanding of how the dynamics were derived, see \cite{williams2022practically}. We introduce the algorithm:
\begin{align}
\begin{split}\label{eqn:th_dyn}
    \dot{\hat \theta} ={}& k \omega_f (-G_J +\\
         & \min \{G_h^{-2},M^+\} \max\{G_J^T G_h - c \eta_h, 0\} G_h )
\end{split}\\
    \dot G_J =& - \omega_f ( G_J - (J( \hat \theta(t) + S(t) ) - \eta_J)M(t) )\label{eqn:gj_dyn}\\
    \dot \eta_J =& -  \omega_f ( \eta_J -J( \hat \theta(t) + S(t) )) \label{eqn:etaj_dyn}\\
    \dot G_h =& -  \omega_f (G_h - (h( \hat \theta(t) + S(t) ) - \eta_h)M(t) )\label{eqn:gh_dyn}\\
    \dot \eta_h =& -  \omega_f ( \eta_h - h( \hat \theta(t) + S(t) )) \label{eqn:etah_dyn} 
\end{align}
where the state variables $\hat \theta, G_J, G_h \in \mathbb{R}^n$, $\eta_J, \eta_h \in \mathbb{R}$. The overall the dimension of the system is $3n+2$. The map is evaluated at $\theta$, defined by 
\begin{equation}
    \theta(t) := \hat \theta(t) + S(t) \; .
\end{equation}
The integer $n$ denotes the number of parameters one wishes to optimize over. The design coefficients are $k, c, \omega_f, M^+\in \mathbb{R}_{>0}$. The perturbation signal $S$ and demodulation signal $M$ are given by  
\begin{align}
    S(t) &= a \left[ \sin(\omega_1 t), \;...\;, \sin(\omega_n t) \right]^T, \\
    M(t) &= \frac{2}{a} \left[  \sin(\omega_1 t), \;...\;, \sin(\omega_n t) \right]^T,
\end{align}
and contain additional design parameters $\omega_i, a \in \mathbb{R}_{>0}$.
\section{Assumptions}
We define
\begin{align}
    \mathcal{C} &= \{\theta \in \mathbb{R}^n : h(\theta) \geq 0 \}, \\
    \partial \mathcal{C} &= \{\theta \in \mathbb{R}^n : h(\theta) = 0 \}, \\
    \mathcal{U}&= \{\theta \in \mathbb{R}^n : h(\theta) \leq 0 \},
\end{align}
where $\mathcal{C}$ is called the `safe set' and $\partial\mathcal{C}$ is its boundary. We also define the notation of a superlevel set of $h$, parameterized by the parameter $\rho \leq 0$ as
\begin{equation}
    \mathcal{C}_\rho = \{\theta \in \mathbb{R}^n : h(\theta) \geq \rho, \rho \in \Ima(h) \cap \mathbb{R}_{\leq 0} \}.
\end{equation}
The sets of the form $\mathcal{C}_\rho$ always contain $\mathcal{C}$ along with some unsafe region given by $\rho$, a non-positive value in the image of $h$. We also use the following assumptions throughout.
\begin{assumption}[Objective Function Conditions] \label{assum:J}
    The objective function $J: \mathbb{R}^n \to \mathbb{R}$ is continuously differentiable with locally Lipschitz Jacobian and satisfies:
    \begin{enumerate}
        \item $\theta^*_c \in \mathcal{C}$ is the unique constrained minimizer of $J$ on $\mathcal{C}$,
        \item if there exists a $\theta$ such that $\nabla J(\theta) = 0$ for $\theta \in \mathcal{C}$, then $\theta=\theta^*_c$.
    \end{enumerate}
\end{assumption}
\begin{assumption}[Barrier Function Conditions] \label{assum:h}
The barrier (or safety) function $h: \mathbb{R}^n \to \mathbb{R}$ is continuously differentiable with locally Lipschitz Jacobian
and satisfies:
\begin{enumerate}    
    \item the safe set $\mathcal{C}$ is non-empty,
    \item for any $\mathcal{C}_\rho$, there exists a $L \in (0, \infty)$ such that $\| \nabla h(\theta)\|>L$ for $\theta \in \mathcal{U} \cap \mathcal{C}_\rho$.\label{assum:h_grad}
\end{enumerate}
\end{assumption}
\begin{assumption}[Optimizer Condition] \label{assum:gradient}
If $ \nabla h ( \theta )^T \nabla J ( \theta ) = || \nabla h ( \theta )|| ||\nabla J ( \theta ) ||$ ($\nabla h ( \theta )$ and $\nabla J (\theta)$ are collinear) for $\theta \in \partial \mathcal{C}$, then $\theta = \theta^*_c$. 
\end{assumption}
\begin{assumption}[Angle Condition] \label{assum:angle}
There exists a $r^*>0$ and  $f^* \in [0,1)$ such that
\begin{equation}
    \frac{\nabla J(\theta)^T \nabla h(\theta)}{|| \nabla J(\theta)|| || \nabla h(\theta)||} \leq f^*,
\end{equation}
for $\theta \in \{\rho \leq h(\theta) \leq 0 \} \cap \{r^* \leq || \theta - \theta^*_c || \}$ for any $\rho \in \Ima(h) \cap \mathbb{R}_{< 0}$.
\end{assumption}
\begin{assumption}[Radial Unboundedness] \label{assum:radially_unbounded}
The function $V = \max \{-h(\theta),0 \} + \max \{J(\theta) - J(\theta^*_c),0\}$ is positive definite and
$ || \theta - \theta^*_c|| \to \infty \implies V \to \infty $.
\end{assumption}
\begin{assumption} [Bounded Levels of $h$] \label{assum:bounded_level}
The family of sets $\mathcal{C}_\rho$ are compact. 
\end{assumption}
\begin{assumption}[ES Constants] \label{assum:es_constants}
The design constants are chosen as $\omega_f, \omega_i, \delta,a, k,c > 0$, where $\omega_i \slash \omega_j$ are rational with frequencies $\omega_i$ chosen such that $\omega_i\neq \omega_j$ and $\omega_i + \omega_j \neq \omega_k$ for distinct $i, j,$ and $k$.
\end{assumption}
\begin{figure}[t]
\centering
\includegraphics[width=1.00\linewidth]{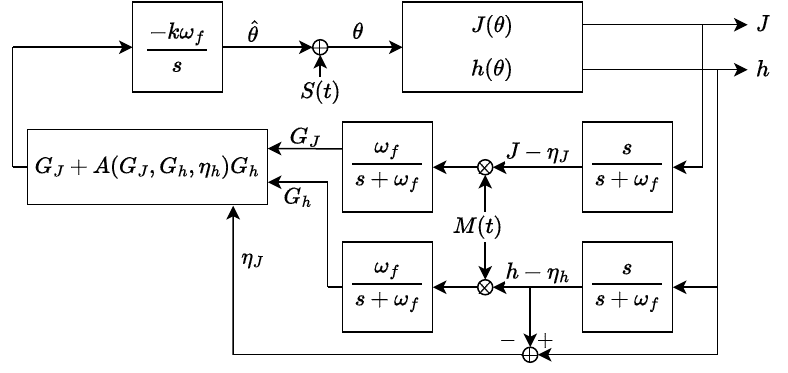}
\caption{Safe Extremum Seeking block diagram. Note that $A(G_J, G_h, \eta_h)=\min \{G_h^{-2},M^+\} \max\{G_J^T G_h - c \eta_h, 0\}$.}
\label{fig:block_diagram}
\end{figure}
Assumptions \ref{assum:J}, \ref{assum:h}, and \ref{assum:gradient} are used to show that the exact dynamics yield only a single equilibrium, and Assumptions \ref{assum:angle}, \ref{assum:radially_unbounded} and \ref{assum:radially_unbounded} are used in the Lyapunov analysis.

Assumption \ref{assum:gradient} is a general condition on the gradients of $J$ and $h$ on the boundary, which is used to force a unique equilibrium of the dynamics. This assumption is congruent with necessary conditions known about solutions in optimization (such as the so called ``method of Lagrange multipliers''). More restrictive assumptions can also be made in place of this. For example, in place of Assumption \ref{assum:gradient} we can assume $\mathcal{C}$ is convex and $J$ is strictly convex on $\mathcal{C}$. Then the dynamics can be shown to yield a unique equilibrium and the convergence analysis proceeds the same otherwise.

Assumption \ref{assum:angle} is concerned with the cosine of the angle between $\nabla J$ and $\nabla h$ far away from the constrained equilibrium in the unsafe set. We utilize this assumption in finding a Lyapunov function with negative time derivative everywhere. It can be shown that this condition is true for $J$ quadratic and convex with $h$ linear - an example of a problem with a semi-infinite safe set. Assumption \ref{assum:radially_unbounded} is used because we use the implication that that the sublevel sets of a particular Lyapunov function (introduced later) can be chosen compact and arbitrarily large.

 We give our main results for either Assumption \ref{assum:angle}-\ref{assum:radially_unbounded} holding or Assumption \ref{assum:bounded_level} holding. This is because Assumption \ref{assum:bounded_level} is strong enough to show SPA stability of the reduced system as it readily yields compact and arbitrarily large invariant sets. If we have a semi-infinite safe set, then we require Assumption \ref{assum:angle}-\ref{assum:radially_unbounded}. 

\section{Global Convergence of the Exact Algorithm}
Before conducting analysis of the ES scheme, we study the optimization algorithm in its' exact form, in order to find the appropriate Lyapunov function which will be used in the next section. Consider the following dynamics:
\begin{multline}
     \dot \theta = F(\theta) = - \nabla J(\theta) + \\ \frac{\nabla h(\theta) }{||\nabla h(\theta) ||^2} \max\{ \nabla J(\theta)^T \nabla h(\theta) - c h(\theta), 0 \}.  \label{eqn:exact_dynamics}
\end{multline} 
We do not risk dividing by zero in the expression of \eqref{eqn:exact_dynamics} as $\nabla h(\theta) = 0 \implies h(\theta)>0$ by Assumption \ref{assum:h}.\ref{assum:h_grad}. Therefore, $\lim_{||\nabla h(\theta)|| \to 0} F(\theta) = -\nabla J(\theta)$ on compact sets.

The differential inequality $ \dot h + ch \geq 0$ is commonly used to show the forward invariance of the safe set \cite{ames2016control}, where it is assumed that the initial condition of any given trajectory is safe. But it also shows attractivity to the safe set. Because when $h(\theta(t)) < 0$, for some $\theta(t)$, then we have $\dot h (\theta(t)) > 0$ and $h$ is increasing in time. So unsafe trajectories become `safer' in an exponential fashion. Therefore, the family of sets $\mathcal{C}_\rho$ is also positively invariant, and not just the case of $\rho = 0$. We state this formally below.
\begin{proposition} \label{prop:forward_inv}
 Under Assumption \ref{assum:h}, the dynamics \eqref{eqn:exact_dynamics} satisfy $\frac{dh(\theta(t))}{dt} + c h(\theta(t)) \geq 0 $ for all $\theta \in \mathbb{R}^n$ and $\mathcal{C}$ is forward invariant. Moreover, all sets $\mathcal{C}_\rho$ are forward invariant and $h(\theta(t)) \geq h(\theta(t_0)) e^{-c t}$ for all $\theta(t_0) \in \mathbb{R}^n$.
\end{proposition}

To prove this, simply compute $\dot h + ch$ and use the fact that $-x + \max \{x,0\} \geq 0$. This proposition will help us establish global convergence of the algorithm for trajectories starting outside of $\mathcal{C}$. The next lemma verifies that all trajectories starting from the safe set, converge to the constrained minimum of $J$.

\begin{lemma} \label{lem:V1}
    Let Assumptions \ref{assum:J}-\ref{assum:gradient} hold. The function $V_1(\theta) = J(\theta) - J(\theta^e)$ is a Lyapunov function for the equilibrium $\theta^e$ on $\mathcal{C}$, yielding strictly $\dot V_1 < 0$ for all $\theta \in \mathcal{C} \setminus \{ \theta^e \}$. The dynamics \eqref{eqn:exact_dynamics} are asymptotically stable for $\theta(t_0) \in \mathcal{C}$.
\end{lemma}

In light of Proposition \ref{prop:forward_inv} (attractivity to $\mathcal{C}$) and Lemma \ref{lem:V1} (convergence within $\mathcal{C}$), it should be intuitive that all trajectories eventually converge on the equilibrium point $\theta^e$. We can demonstrate this fact with the following Lyapunov argument on any compact, invariant set $\mathcal{C}_\rho$ around the equilibrium. A key idea used is that any initial condition is an element of some set of the form $\mathcal{C}_\rho$.

\begin{lemma} \label{lem:global_conv_exact}
Let Assumptions \ref{assum:J}-\ref{assum:gradient} hold, and let either Assumptions \ref{assum:angle}-\ref{assum:radially_unbounded} or Assumption \ref{assum:bounded_level}. Consider the Lypaunov function
\begin{equation}
    V(\theta) =  \max\{-\alpha h(\theta), 0\} + \max \{J(\theta) - J(\theta^e),0\}.
\end{equation}
For any $\mathcal{C}_\rho$, there exists $\alpha \in (0, \infty)$ such that $\dot V(\theta) < 0$ for $\theta \in \mathcal{C}_\rho \setminus \{\theta^e\}$. The dynamics \eqref{eqn:exact_dynamics} are globally asymptotically stable.
\end{lemma}

To prove the result in Lemma \ref{lem:global_conv_exact}, time derivatives of the Lyapunov function $V$ must be computed in all regions of the state space within any arbitrary invariant set $\mathcal{C}_\rho$. In this paper we will sketch part of the proof and compute the time derivative of $V$ in the case where $\theta \in \mathcal{U} \cap \{J(\theta) - J(\theta^e) \geq 0\} \cap \{\nabla J(\theta)^T \nabla h(\theta) - c h(\theta) \geq 0\}$  (the $\max$ term is active in the dynamics).

First compute expressions for $\dot V$:
\begin{align}
    \dot V =& -\alpha c |h(\theta)| + \dot{V}_1, \nonumber \\
\begin{split} \label{eqn:case4_Vdot}
    ={}& c |h(\theta)| \left( -\alpha  + \frac{\nabla J(\theta)^T \nabla h(\theta)}{\|\nabla h(\theta) \|^2} \right)  - \\ & \nabla J(\theta)^T \left(I - \frac{\nabla h(\theta) \nabla h(\theta)^T}{\|\nabla h(\theta) \|^2} \right) \nabla J(\theta).
\end{split}
\end{align}
Denoting 
\begin{equation}
f(\theta) := \nabla J(\theta)^T \nabla h(\theta) /(||\nabla J(\theta) || || \nabla h(\theta)||) \leq 1 ,   
\end{equation}
we can rewrite and bound $\dot V$ as
\begin{align}
\begin{split}
    \dot V ={}& -\alpha c |h(\theta)| - (1-f^2(\theta)) || \nabla J(\theta) ||^2 + \\ &f(\theta) \frac{c|h(\theta)|}{||\nabla h(\theta) ||} ||\nabla J(\theta) ||,
\end{split}  \nonumber\\
 \leq &  c |h(\theta)| \left( -\alpha + \frac{||\nabla J(\theta) ||}{L} \right) - (1-f^2(\theta)) || \nabla J(\theta) ||^2. \label{eqn:form_1} 
\end{align}
The existence of $L$ is assumed by Assumption \ref{assum:h} and recall we assume that $\theta \in \mathcal{C}_\rho$ for some $\mathcal{C}_\rho$.

\textbf{Case A) Assumption \ref{assum:bounded_level}:} Consider the case of Assumption \ref{assum:bounded_level}, then $\mathcal{C}_\rho$ is compact and we can choose
\begin{equation}
    \alpha > L^{-1} \sup_{\theta \in \mathcal{C}_\rho} \| \nabla J(\theta) \| \label{eqn:alpha_C_compact},
\end{equation}
which yields $\dot V \leq 0$ in \eqref{eqn:form_1}.

\textbf{Case B) Assumptions \ref{assum:angle} - \ref{assum:radially_unbounded}:} Consider the case of Assumptions \ref{assum:angle} - \ref{assum:radially_unbounded}. If $||\nabla J(\theta)||$ is bounded on $\mathcal{C}_\rho$ then we can choose $\alpha$ as in \eqref{eqn:alpha_C_compact}. So consider the nontrivial case where $|| \nabla J(\theta)||$ unbounded on $\mathcal{C}_\rho$. Then from Assumption \ref{assum:angle}, there exists a scalar $f^*$ and a compact set $\Omega$ containing $\theta^e$ such that $0<f(\theta) \leq f^* < 1$ for all $\theta \notin \Omega$. Therefore letting $\tilde{f} =  1-{f^*}^{2}$ with $\tilde{f} \in (0,1)$ we have
\begin{align}
    \dot V & \leq  c |h(\theta)|\left(-\alpha + \frac{||\nabla J(\theta) ||}{L}\right) - \tilde f || \nabla J(\theta) ||^2 
\end{align}
for $\theta \notin \Omega$, which also yields the bound
\begin{align}
    \dot V & \leq  -\alpha c |h(\theta)| + \frac{c |\rho| }{L} ||\nabla J(\theta) || - \tilde f || \nabla J(\theta) ||^2.
\end{align}
This implies when $||\nabla J(\theta) ||>\frac{|\rho| c}{L \tilde f}$, then $\dot V < 0$ for $\theta \notin \Omega$. Therefore choosing
\begin{equation}
    \alpha > \max \left\{ L^{-1} \sup_{\theta \in \Omega} \| \nabla J(\theta) \|, \frac{c |\rho|}{L^2 \tilde f} \right\}
\end{equation}
yields $\dot V \leq 0$.

For a complete proof, which lies outside the scope of this paper, one must also consider computing the time derivative of $V$ in the other regions of the state space. 

This Lyapunov function has an interesting connection to literature. We note that a closely related Lyapunov function for `gradient flow' systems in a more general setting was also discovered and can also be used to show stability \cite{allibhoy2022control}. 

\section{Safe Extremum Seeking of Static Maps}
The basic outline of this section starts by following the framework in \cite{nesic2010unifying}. As a preliminary, we first present Lemma \ref{lem:additive_disturbance} to aid the reader in understanding later results. Then, we perform a series of transformations starting from the original system, resulting in the construction of a `reduced model'. We then use a Lyapunov argument and the help of Lemma \ref{lem:global_conv_exact} to show the reduced model is SPA stable, which yields the original system SPA stable. Finally, we present our notion of practical safety which follows from the SPA stability of the ES scheme. 
\begin{remark}
Authors in \cite{nesic2010unifying} require that the reduced model be ``robust'' to ``disturbances'' (SPA stable) which arise as a result of the averaging procedure. We use the terminology ``robust'' and ``disturbance'' in the following Lemma and the current section to refer to this fictitious disturbance of the reduced system. This reduced system is fully constructed later in \eqref{eqn:reduced} and is required for the analysis.
\end{remark}
Consider the estimated quantities:
\begin{equation}
    Q(\theta)^T := [\nabla J (\theta)^T, \nabla h (\theta)^T, h(\theta)].
\end{equation}
The dynamics in \eqref{eqn:exact_dynamics} can be thought of as a function of the estimated variables $F = F(Q(\theta))$. And with a small disturbance $w(t)^T = [w_1(t)^T, w_2(t)^T, w_3(t)]$ we can consider a disturbed system written as 
\begin{equation}
    \dot \theta = F(Q(\theta) + w(t)) . \label{eqn:disturbed_dyn}
\end{equation}
The next result is a useful preliminary as it says that the algorithm, without perfect estimates of the gradients, can be written as the exact version of the algorithm with a small additive disturbance. It is also useful as it shows that this holds for a large enough $M^+$.
\begin{lemma} [Additive Disturbance]\label{lem:additive_disturbance}
        Under Assumptions \ref{assum:J} - \ref{assum:h}, for any compact set $\Omega$ there exists a $M^+, \epsilon^*>0$ such that for all $\epsilon \in (0, \epsilon^*) $ and $||w || < \epsilon$, the disturbed dynamics in \eqref{eqn:disturbed_dyn} can be written as
    \begin{multline}
        F(Q(\theta) + w(t)) = - \nabla J(\theta) + \\ \frac{\nabla h(\theta) }{||\nabla h(\theta) ||^2} \max\{ \nabla J(\theta)^T \nabla h(\theta) - c h(\theta), 0 \} + O(\epsilon).
    \end{multline}
for all  $\theta \in \Omega$.
\end{lemma}
The proof of Lemma \ref{lem:additive_disturbance} is outside the scope of the paper. We now turn our attention to the original dynamics in \eqref{eqn:th_dyn} - \eqref{eqn:etah_dyn}, and make a series of transformations, following the general ideas in \cite{nesic2010unifying}. Defining
\begin{align}
    F_0\left(\xi\right) &:= -\xi_1 + \min \{\xi_2^{-2},M^+\} \max\{\xi_1^T \xi_2 - c \xi_4, 0\} \xi_3,
\end{align}
with
\begin{align}
    \xi^T &:= [G_J^T, \eta_J, G_h^T, \eta_h],\\
    \zeta^T &:= [(J(\theta) - \xi_2)M(t)^T, J(\theta), (h(\theta) - \xi_4)M(t)^T, h(\theta)],
\end{align}
we can rewrite \eqref{eqn:th_dyn} - \eqref{eqn:etah_dyn} as
\begin{align}
    \dot{\hat{\theta}} &= k \omega_f F_0(\xi), \label{eqn:theta_dyn_rewrite}\\
    \dot \xi &= -\omega_f (\xi - \zeta(t,\theta, \xi, a)),  \label{eqn:estimate_dyn_rewrite}
\end{align}
recalling $\theta = \hat \theta + S(t)$. Letting $\tilde \theta = \hat \theta - \theta^*_c$ and $\tau = \omega_f t$, the system in the new time scale is
\begin{align}
    \frac{d \tilde \theta}{d \tau} &= k F_0(\xi), \label{eqn:theta_dyn_tau}\\
    \frac{d \xi}{d \tau} &= -\left(\xi - \zeta \left( \frac{\tau}{\omega_f}, \tilde \theta + \theta^*_c + S\left(\frac{\tau}{\omega_f}\right), \xi, a\right) \right).\label{eqn:xi_dyn_tau}
\end{align}
We can take the average of the system (see \cite{Khalil}, \cite{nesic2010unifying}) to compute
\begin{align}
     \frac{d \tilde \theta_{av}}{d \tau} =& k F_0(\xi_{av}), \label{eqn:theta_dyn_tau_av}\\
    \frac{d \xi_{av}}{d \tau} =& -\left(\xi_{av} - \mu(\tilde \theta_{av},a) \right), \label{eqn:xi_dyn_tau_av}\\
    \begin{split} \label{eqn:D_def}
        D(\tilde \theta_{av}+ \theta^*_c) :={}& [\nabla J(\tilde \theta_{av}+ \theta^*_c)^T, J(\tilde \theta_{av}+ \theta^*_c), \\ & \nabla h(\tilde \theta_{av}+ \theta^*_c)^T, h(\tilde \theta_{av}+ \theta^*_c)]^T ,  
    \end{split} \\    
    \mu(\tilde \theta_{av},a) :=& D(\tilde \theta_{av}+ \theta^*_c) + O(a). \label{eqn:mu_def}
\end{align}
Making another time transformation $s = k \tau$ we have
\begin{align}
     \frac{d \tilde \theta_{av}}{d s} &= F_0(\xi_{av}), \label{eqn:theta_dyn_s}\\
    k \frac{d \xi_{av}}{d s} &= -\left(\xi_{av} - \mu(\tilde \theta_{av},a) \right) .\label{eqn:xi_dyn_s}
\end{align}
Taking $k=0$, we can derive the singularly perturbed (or reduced) system with a quasi steady state
\begin{equation}
    z_s := \xi_{av} = \mu(\tilde \theta_{av},a).
\end{equation}
Defining
\begin{equation}
    y = \xi_{av} - \mu(\tilde \theta_{av},a),
\end{equation}
the boundary layer system (with $\tau = s/k$) is
\begin{equation}
    \frac{dy}{d \tau} = -\left(\xi_{av} - \mu(\tilde \theta_{av},a) \right) = -y.
\end{equation}
The boundary layer system is UGAS uniformly in $\xi_{av}$ and $t_0$. The reduced system is
\begin{equation}
    \frac{d \theta_{r}}{d s} = F_0(\mu(\theta_{r},a)) = F_0(D(\theta_{r}+ \theta^*_c) + O(a)). \label{eqn:reduced}
\end{equation}

The reduced system \eqref{eqn:reduced} is SPA stable in $a$. This is due to the following argument: take any positive pair $(\Delta, \nu)$, and $\theta_r(0) \leq \Delta$. Using Lemma \ref{lem:global_conv_exact} and \ref{lem:additive_disturbance} we can find a Lyapunov function $V$ such that 
\begin{align}
    & \dot V(\theta_r) = -W(\theta_r) + O(a), \\
    & \dot h(\theta_r) + ch(\theta_r) \geq O(a) ,
\end{align}
for strictly positive definite $V$ and $W$, for a sufficiently large $M^+$ on any compact set. In Lemma \ref{lem:global_conv_exact} we considered two cases: 1) Assumption \ref{assum:bounded_level} (all $\mathcal{C}_\rho$ are compact) 2) Assumptions \ref{assum:angle}-\ref{assum:radially_unbounded} ($V$ is radially unbounded and the angle condition is obeyed). In either case we can construct $\mathcal{S}$, a compact invariant set containing $B(0)_\Delta$.

 Either $\mathcal{S} = \mathcal{C}_\rho$, for some $\mathcal{C}_\rho$ in case 1 or $\mathcal{S} = \mathcal{C}_\rho \cap \{V(\theta_r) \leq \bar V \}$, for some $\bar V>0$ in case 2 for sufficiently small $a$ and sufficiently large $M^+$. One may then choose an $a$ even smaller to achieve $\lim_{t \to \infty} ||\theta_r(t)|| < \nu$ for any $\nu>0$. Therefore there exists an $a^*$ such that for any $a \in (0, a^*)$ the solutions of \eqref{eqn:reduced} satisfy
\begin{equation}
    ||\theta_r(s)|| \leq \beta_{\theta}(|| \theta_r(0) ||, s) + \nu
\end{equation}
for some $\beta_{\theta} \in \mathcal{KL}$, $M^+>0$ and for all $||\theta_r(0)|| \leq \Delta$. SPA stability of \eqref{eqn:reduced} satisfies Assumption 2 in \cite{nesic2010unifying}).

Using \cite[Theorem 1]{nesic2010unifying} we make the conclusion below in Theorem \ref{thm:spa_stable}. Note: The authors in \cite{nesic2010unifying} prove this result based upon the earlier work \cite[Lemma 1,2]{tan2005non}. \cite[Lemma 2]{tan2005non} concludes the SPA stability in $[a, k]$ of \eqref{eqn:theta_dyn_s}-\eqref{eqn:xi_dyn_s} from the SPA stability in $a$ of \eqref{eqn:reduced}. \cite[Lemma 1]{tan2005non} concludes the SPA stability of \eqref{eqn:theta_dyn_rewrite}-\eqref{eqn:estimate_dyn_rewrite} from SPA stability in $[a, k]$ of \eqref{eqn:theta_dyn_s}-\eqref{eqn:xi_dyn_s}. 
Let
\begin{equation}
    z = \xi - \mu(\tilde \theta, a).
\end{equation}

\begin{theorem} [Semi-Global Practical Stability] \label{thm:spa_stable}
Let Assumptions \ref{assum:J}-\ref{assum:gradient}, \ref{assum:es_constants} hold. Also, let either \ref{assum:angle} and \ref{assum:radially_unbounded} hold or \ref{assum:bounded_level} hold. Then there exists $\beta_\theta, \beta_\xi \in \mathcal{KL}$ such that: for any positive pair $(\Delta, \nu)$ there exist $M^+, \omega_f^*, a^*>0$, such that for any $\omega_f \in (0, \omega_f^*)$, $a \in (0,a^*)$, there exists $k^*(a)>0$ such that for any $k \in (0, k^*(a))$ the solutions to \eqref{eqn:theta_dyn_rewrite}-\eqref{eqn:estimate_dyn_rewrite} satisfy
\begin{align}
    || \tilde \theta(t) || &\leq \beta_\theta \left( ||\tilde \theta (t_0)||, k \cdot \omega_f \cdot (t-t_0) \right) + \nu, \label{eqn:theta_KL_bound} \\
    || z(t) || &\leq \beta_\xi \left( ||z (t_0)||, \omega_f \cdot (t-t_0) \right)  + \nu, \label{eqn:z_KL_bound}
\end{align}
for all $||[\tilde \theta(t_0)^T, z(t_0)^T ]^T|| \leq \Delta$, and all $t \geq t_0\geq 0$.
\end{theorem}

Equation \eqref{eqn:z_KL_bound} tell us about convergence of estimated quantities $\xi(t)$ to their exact values $ D(\tilde \theta(t)  + \theta^*_c)$. Using  \eqref{eqn:mu_def} we can write
\begin{align}
    e &:= \xi - D(\tilde \theta  + \theta^*_c), \label{eqn:e_def}\\
     || e(t) || &\leq \beta_\xi (||z(t_0)||, \omega_f(t-t_0)) + v + |O(a)|. \label{eqn:e_bound}
\end{align}
The variable $e$ can be though of as the estimator error of the various measurements and gradients of the static maps. The bound above says that the estimated quantities can be made to converge arbitrarily close to their true values, at a time scale faster than the movement of the parameter $\tilde \theta$ (which is governed by the gain $k$). Note, the functions $\beta_\theta, \beta_\xi$ are independent of $a, k, \omega_f$ \cite{tan2005non}. Because $\dot{\tilde \theta}$ is proportional to $k$ \eqref{eqn:theta_dyn_rewrite}, we can choose a $k$ such that the change in $\tilde \theta (t)$ over some time interval can be small relative to the change in $e(t)$ over the same interval. From Lemma \ref{lem:additive_disturbance} we know that there must be a small enough disturbance in the estimated quantities such that the dynamics can be written linearly, but only after the transient of $\beta_\xi$ has been sufficiently diminished. Therefore, we make the following claim.
\begin{theorem}[Semi-Global Practical Safety] \label{thm:prac_safety}
    Suppose Theorem \ref{thm:spa_stable} holds. For any $\Delta>0$ there exists $\delta^*>0$ such that for any $\delta \in (0, \delta^*)$ : there exists $M^+, a^{**}, \omega_f^{**}>0$ such that for any $a \in (0, a^{**})$, $\omega_f \in (0, \omega_f^{**})$, there exists $k^{**}(a)>0$ such that for any $k \in (0, k^{**}(a))$,
\begin{equation}
    h(\theta(t))  \geq h(\theta(t_0)) e^{-c k \omega_f (t-t_0)} + O(\delta) \label{eqn:practical_safety_inequality}
\end{equation}
     for all $||[\tilde \theta(t_0)^T, z(t_0)^T ]^T|| \leq \Delta$ for all $t \in [t_0, \infty]$. 
\end{theorem}
Idea of proof:
Consider the system in the form of \eqref{eqn:theta_dyn_rewrite}-\eqref{eqn:estimate_dyn_rewrite}, and also recall that $z=\xi - \mu(\tilde \theta, a)$. From Theorem \ref{thm:spa_stable} it was shown that the solutions satisfy the bounds in \eqref{eqn:theta_KL_bound} - \eqref{eqn:z_KL_bound}, with rates of decay $k \omega_f$ and $\omega_f$ respectively. We choose $M^+, \omega_f^*, a^*>0$ such that the $\nu$ is sufficiently small, and the estimator state $\xi$ converges to a region which is close to $\mu(\tilde \theta, a)$. Next, one can further restrict $a$ such that $\mu(\tilde \theta, a)$ is close to true gradients. Therefore, one can argue that after some finite time $T$, $\xi(t)$ will be within some small error away from the true gradients. Therefore the system, after this transient, \eqref{eqn:theta_dyn_rewrite} can be written as $\dot{\hat{\theta}} = k \omega_f F_0 (D(\hat{\theta}) + e(t))$. From Lemma \ref{lem:additive_disturbance}, we can write $F_0 (D(\hat{\theta}) + e(t)) = F_0 (D(\hat{\theta})) + O(\delta)$ for $||e(t)|| \leq \delta$. Therefore for t in $[T, \infty)$ we have $\dot h + c h \geq O(\delta)$, implying \eqref{eqn:practical_safety_inequality}. During the finite time from $[0, T]$ we can bound changes in $\hat{\theta}$ with $k$ as the time $T$ is independent of $k$, implying \eqref{eqn:practical_safety_inequality}. 
Note this argument makes use of the crucial fact that the functions $\beta_\theta, \beta_\xi$ are independent of $a, k, \omega_f$ \cite{tan2005non}.

This argument implies the intervals of $a, \omega_f, k$ given in Theorem \ref{thm:prac_safety} are perhaps a more strict set of intervals than the ones given in Theorem \ref{thm:spa_stable} - if the user desires the type of safety given in \eqref{eqn:practical_safety_inequality}. Nonetheless, the safety result is elegantly analogous to the statement on stability. Theorem \ref{thm:spa_stable} says that for any set of initial conditions, one should be able to adjust $a, \omega_f, k$ such that trajectories are $\nu$-practically stable. Theorem \ref{thm:prac_safety} says that for any set of initial conditions, one should be able to adjust $a, \omega_f, k$ such that trajectories are $\delta$-practically stable.

\begin{figure}[t!]
\centering
\includegraphics[width=0.95\linewidth]{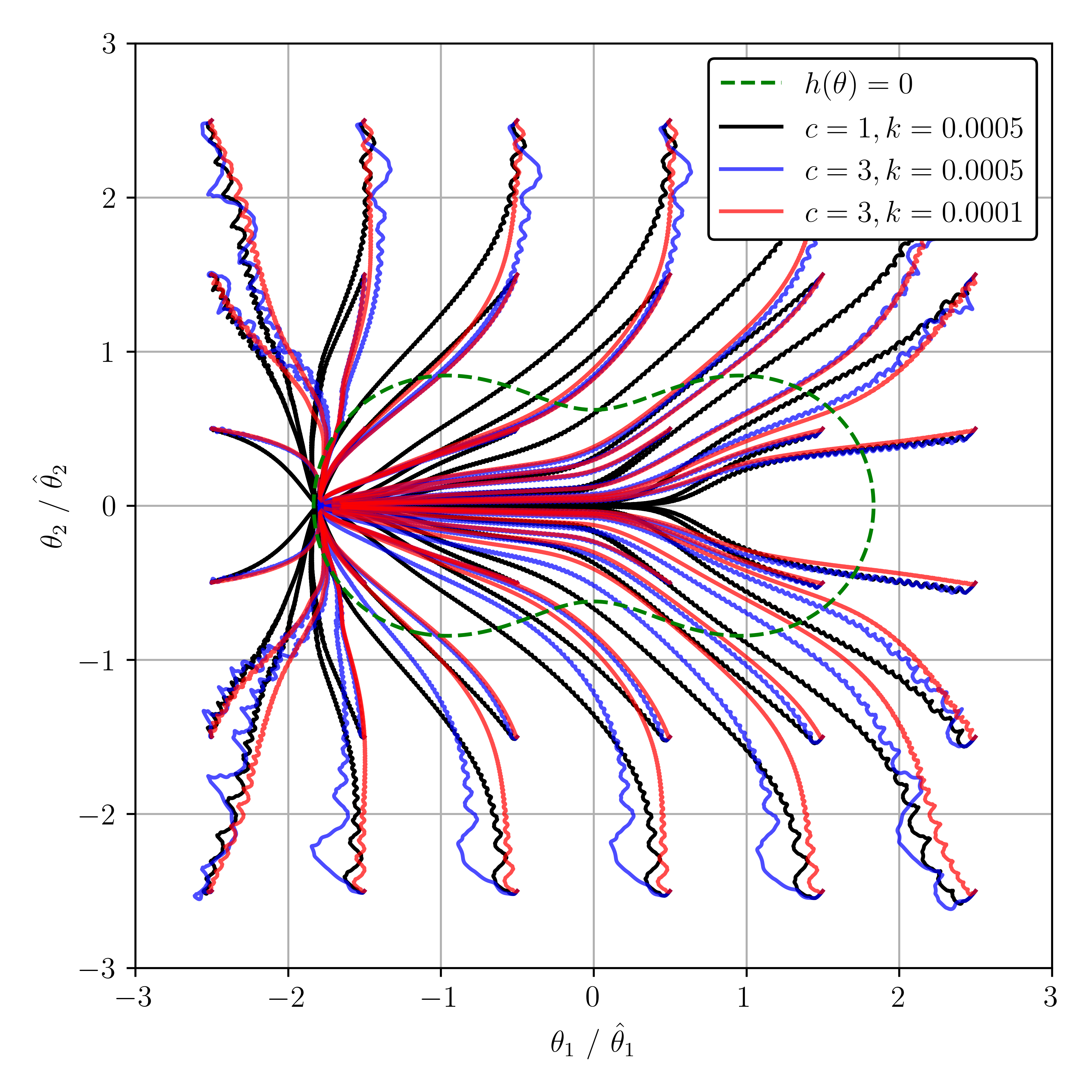}
\caption{Simulation of the algorithm with various $c$ and $k$.}
\label{fig:ex1}
\end{figure}

\section{Example}
Consider the following maps and parameters: $J(\theta) = (\theta_1 + 3)^2 + \theta_2^2$, $h(\theta) = e^{-(\theta_1-1)^2 - \theta_2^2} + e^{-(\theta_1+1)^2 - \theta_2^2} -0.5$, $a= 0.1 $, $\omega_f= 10$, $M^+= 10,000$, $\omega_1= 10$, $\omega_2= 13$. The safety function chosen has the characteristic that all the level sets for $h\leq0$ ($\mathcal{C}_\rho$) are compact, and there is a unique minimum of $J$ on $\mathcal{C}$. The levels of $J$ and $h$ can also be shown to obey the optimizer criterion in Assumption \ref{assum:gradient}. Note that the safe set is non convex. In Fig. \ref{fig:ex1} we simulate a grid of initial conditions (all estimator states $G_J, \eta_h, G_h, \eta_h$ are initialized to zero) with three scenarios a) $c=1$, $k=0.0005$ b) $c=3$, $k=0.0005$ c) $c=3$, $k=0.0001$. 

In Fig. \ref{fig:ex1}, notice that trajectories converge to a region around the minimizer of $J$ on $\mathcal{C}$. Also, trajectories that enter the safe set always remain in the safe set. Consider the effect of $c$. When in the safe set, a smaller $c$ restricts the approach towards the barrier. This can be see with black trajectories favoring the center-line (the safest part of the set) over the blue/red around the region $\hat \theta \approx [0.5, 0.0]$. Outside the safe set, compare black and blue/red trajectories at $\hat \theta(0) = [-2.5, -0.5]$. Here, a higher constant $c$ dictates a slow escape towards the boundary of the safe set which is shown by the blue/red trajectories going first towards the barrier, before converging. Therefore, higher $c$ has the effect that the optimization of $J$ is favored in the safe region, but in the unsafe region, safety is favored. 

Note the increased transient wiggles introduced at the start of the blue trajectories (compared to black), when increasing $c$ with no adjustment in $k$. We can fix this sign of instability by lowering $k$, and achieve the more smooth trajectories in red. This may be necessary because changing $c$ fundamentally changes the dynamics, and the same gain $k$ may no longer be appropriate for SPA stability - although in this example the stability still remains for the region we have chosen. 
 
\section{Conclusion}
 We have introduced an ES algorithm that not only minimizes an objective but does so with a guarantee of practical safety. A Lyapunov analysis shows that for a semi-global region of attraction there exist intervals of design coefficients such that both practical safety and practical stability exist. Our example demonstrates a convergence of a non-convex problem.

\bibliographystyle{IEEEtranS}
\bibliography{refs}

\end{document}